\documentclass[11pt]{article}
\usepackage{amsmath}

\newcommand{\R}{{\sf I\hspace{-.15em}R}}

\def\di{\displaystyle}

\def\bs{\bigskip}

\def\saa{\Sigma_A^+}
\def\sa{\Sigma_A}

\def\ff{{\mathcal F}}

\def\Pr{\mbox{\rm Pr}}
\def\ms{\medskip}
\def\var{\mbox{\rm var}}
\def\be{\begin{equation}}
\def\ee{\end{equation}}
\def\beqn{\begin{eqnarray}}
\def\eeqn{\end{eqnarray}}
\def\endofproof{{\rule{6pt}{6pt}}}
\def\MM{{\mathcal M}}

\def\tpsi{\tilde{\psi}}
\def\hf{\hat{f}}
\def\Int{\mbox{\rm Int}}
\def\II{{\mathcal I}}

\begin{document}

\noindent
{\bf\large A uniform estimate for rate functions in large deviations}

\bs

\noindent
{\bf Luchezar Stoyanov}\footnote{
{\footnotesize School of Mathematics and Statistics, University of Western Australia, 35 Stirling Hwy, Crawley 6009 WA, Perth, Australia;
E-mail: luchezar.stoyanov@uwa.edu.au}}

\vspace{2cm}

\footnotesize

\noindent
{\bf Abstract.} Given  H\"older continuous functions $f$ and $\psi$ on a sub-shift of finite type $\saa$ 
such that $\psi$ is not cohomologous to a constant, the classical large deviation principle holds (\cite{OP}, \cite{Kif}, \cite{Y}) 
with a rate function $I_\psi\geq 0$ such that $I_\psi (p) = 0$ iff $p =  \int \psi \, d \mu$, where $\mu = \mu_f$ is the equilibrium state of $f$.
 In this paper we derive a uniform estimate from below for $I_\psi$ 
 for $p$ outside an interval containing $\tpsi = \int \psi \, d\mu$, which depends 
only on the sub-shift $\saa$, the function $f$, the norm $|\psi|_\infty$, the H\"older constant of $\psi$ and
the integral $\tpsi$. Similar results can be derived in the same way e.g. for Axiom A diffeomorphisms on basic sets.

\bs


\bs


\normalsize

\section{Introduction} \label{sect1}
\renewcommand{\theequation}{\arabic{section}.\arabic{equation}}

Let  $T: X \longrightarrow X$ be a transformation preserving an ergodic probability measure $\mu$ on a set $X$.
Given an observable $\psi: X \longrightarrow \R$, Birkhoff's ergodic theorem implies that
$$\frac{\psi_n(x)}{n} = \frac{\psi(x) + \psi(T(x)) + \ldots + \psi(T^{n-1}(x))}{n} \to \int_X \psi \, d\mu$$
for $\mu$-almost all $x\in X$. It follows from general large deviation principles (see \cite{OP}, \cite{Kif},
\cite{Y}) that if $X$ is a mixing basic set for an Axiom A diffeomorphism $T$, and
$f$ and $\psi$ are H\"older continuous functions on $X$ with {\it equilibrium states}
$\mu = \mu_f$ and $\mu_\psi$, respectively, and $\psi$ is not cohomologous to a constant (see the definition below), 
then there exists  a real-analytic {\it rate function} $I = I_\psi: \Int(\II_\psi) \longrightarrow [0,\infty)$,
where $\di \II_\psi = \left\{ \int \psi\, d m : m \in \MM_T \right\}$, such that
\be \label{eq:1.1}
\lim_{\delta \to 0} \lim_{n \to \infty} \frac{1}{n} \log \mu \left( \left\{ x\in X : \frac{\psi_n(x)}{n}
\in (p - \delta , p + \delta)\right\}\right) = - I_\psi (p)
\ee
for all $ p \in \Int ({\mathcal I}_{\psi})$. Here $\MM_T$ is the {\it set of all $T$-invariant Borel probability measures} on $X$.
Moreover, $I(p) = 0$ if and only if $p = \int \psi \, d \mu$, and the (closed) interval ${\mathcal I}_{\psi}$ is non trivial,
since $\psi$ is not cohomologous to a constant. 

Similar large deviation principles apply for any  subshift of finite type  $\sigma : \saa \longrightarrow \saa$ on
an one-sided shift space
$$\Sigma_A^+ = \{ \xi = (\xi_0, \xi_1, \ldots,\xi_m,\ldots) : 1\leq \xi_i \leq s_0 , A(\xi_i,\xi_{i+1}) = 1
\; \mbox{\rm for all } \; i\geq 0\;\} .$$
Here $A$ is a $s_0\times s_0$ matrix of $0$'s and $1$'s ($s_0\geq 2$). We assume that $A$ is 
{\it aperiodic}, i.e.  there exists an integer $M > 0$ such that  $A^M(i,j) > 0$ for all $i$, $j$ (see
e.g. Ch. 1 in \cite{PP}). The {\it shift map} $\sigma$ is defined by $\sigma(\xi) = \xi'$, 
where $\xi'_i = \xi_{i+1}$ for all $i \geq 0$. We consider $\saa$ with a metric $d_\theta$ defined for
some constant $\theta \in (0,1)$ by $d_\theta(\xi,\eta) = 0$ if $\xi = \eta$ and $d_\theta(\xi,\eta) = \theta^k$ if 
$\xi \neq \eta$ and $k\geq 0$ is the maximal integer with $\xi_i = \eta_i$ for $0\leq i \leq k$.

For any function  $g: \saa \longrightarrow \R$ set
$$\var_k g = \sup\{ |g(\xi)- g(\eta)| : \xi_i = \eta_i ,\: 0\leq i \leq k\}\:\:\:, \:\:\:
|g|_\theta = \sup \left\{ \frac{\var_k g}{\theta^k} : k \geq 0\right\} \:,$$
$$|g|_\infty = \sup \{ |g(\xi)| : \xi\in\saa\}\:\:\:, \:\:\: \| g\|_\theta = |g|_\theta + |g|_\infty\;.$$ 
Denote by $\ff_\theta(\saa)$ the {\it space of all functions} $g$ on $\saa$ with $\|g\|_\theta < \infty$.

Two functions $f,g $ on $\ff_\theta(\saa)$ are called {\it cohomologous} if there exists a continuous
function $h$ on $\saa$ such that $f = g + h\circ \sigma - h$.

The {\it Ruelle transfer operator} $L_f : C(\saa) \longrightarrow C(\saa)$  is defined by 
$$L_fg(x) = \sum_{\sigma(y) = x} e^{f(y)}\, g(y)\;.$$ 
Here $C(\saa)$ denotes the space of all continuous functions $g : \saa \longrightarrow \R$ with respect to the metric $d_\theta$.
Denote by $\Pr(\psi)$ the {\it topological pressure} 
$$\Pr (\psi) = \sup_{m \in {\mathcal M_{\sigma}}} \Bigl( h_{\sigma}(m) + \int \psi \; d m \Bigr) $$
of $\psi$ with respect to the map $\sigma$, where ${\mathcal M}_{\sigma}$ is the {\it set of all $\sigma$-invariant 
probability measures on $\saa$ } and $h_{\sigma}(m)$ is the {\it measure theoretic entropy} of $m$ with respect to 
$\sigma$ (see \cite{PP} or \cite{R}). Given $\psi\in \ff_\theta(\saa)$, there exists a unique $\sigma$-invariant probability measure $\mu_\psi$ on
$\saa$ such that 
$$\Pr (\psi) = h_{\sigma}(\mu_\psi) + \int \psi \; d \mu_\psi $$
(see e.g. Theorem 3.5 in \cite{PP}).
The measure $\mu_\psi$ is called the {\it equilibrium state} of $\psi$.

For brevity throughout  we write $\di \int h\, d m$ for $\di \int_{\saa} h\, dm $.

In what follows we assume that $\theta \in (0,1)$ is {\bf fixed constant, $f : \saa \longrightarrow \R$
is a  fixed function in $\ff_\theta(\saa)$ and $\mu = \mu_{f}$.}

As we mentioned earlier, it follows from the Large Deviation Theorem (\cite{Kif}, \cite{OP}, \cite{Y}) that if $\psi$ is not 
cohomologous to a constant, then there exists a real analytic {\it rate function}
$I = I_\psi: \Int(\II_\psi) \longrightarrow [0,\infty)$ with $I(p) = 0$ iff
$p = \int \psi\, d \mu$ for which (\ref{eq:1.1}) holds. More precisely, we have
\be \label{eq:1.2}
- I(p) = \inf \{\Pr (f + q \psi) - \Pr(f) - q p: \: q \in \R\} .
\ee
It is also known that
\be \label{eq:1.3}
\left[\frac{d }{dq}\, \Pr(f + q\psi)\right]_{q= \eta} = \int \psi \, d\mu_{f +\eta \psi} ,
\ee
and $\Pr(f + q \psi)$ is a strictly convex function of $q$ (see \cite{PP}, \cite{R} or \cite{La}).

In his paper we derive an estimate from below for $I_\psi(p)$ for $p$ outside an interval containing 
$$\tpsi = \int \psi \, d\mu .$$
The estimate depends only on $|\psi|_\infty$, $\tpsi$, $|\psi|_\theta$ and some constants determined by the
given function $f$. 
In what follows we use the notation $\min \psi = \min_{x\in \saa} \psi(x)$,
$$b = b_\psi = \max \{ 1, |\psi|_\theta \} \quad, \quad B_\psi = \tpsi - \min \{ 0, \min\psi\} .$$
Since $\tpsi > \min \psi$ ($\psi$ is not cohomologous to a constant),
we have $\tpsi - \min\psi > 0$, so $B_\psi > 0$ always.

\bs

\noindent
{\bf Theorem 1.} {\it  Let $f, \psi \in \ff_\theta(\saa)$ be real-valued functions. Assume that $\psi$ is not cohomologous
to a constant, and let $0 < \delta_0 <  B_\psi$.  Then for all $p \notin [\tpsi - \delta_0, \tpsi + \delta_0]$ we have 
$$I_\psi (p) \geq  \frac{\delta_0 q_0}{2} ,$$
where $q_0 =  \min \left\{ C , \frac{1}{b} \right\}$ for some constant $C >0$ depending only on $|f- \Pr(f)|_\infty$, 
$|f|_\theta$,  $|\psi|_\infty$,  $\tpsi$ and $\delta_0$. }

\bs

The motivation to try to obtain estimates of the kind presented in Theorem 1 comes from attempts to get some kind of an 
`approximate  large deviation principle' for characteristic functions $\chi_K$ of arbitrary compact sets $K$ of positive measure. 
In the special case when the boundary $\partial K$ of $K$ is `relatively regular' (e.g. $\mu(\partial K) = 0$) large deviation results 
were established by Leplaideur and Saussol in \cite{LS}, and also by Kachurovskii and Podvigin   \cite{KP}. 
The next example presents a first step in the case of an  arbitrary compact  set $K$ of positive measure.

\ms

\noindent
{\bf Example 1.} Let $K$ be a compact subset of $\saa$ with $0 < \mu(K) < 1$, let $0 < \delta_0 \leq \mu(K)$, and let $\psi$ be a 
H\"older continuous function that approximates $\chi_K$ from above, i.e. $0 \leq \psi \leq 1$, $\psi = 1$ on $K$ and 
$\psi = 0$ outside 
a small neighbourhood $V$ of $K$. Then $b = |\psi|_\theta > > 1$ if $V$ is sufficiently small, so $q_0$ in Theorem 1 has the form 
$q_0 = 1/b$. It then follows from Theorem 1 (in fact, from Lemma 1 below) that $I_\psi (p) \geq \frac{\delta_0}{2 |\psi|_\theta}$ for 
$p \notin [\tpsi-\delta_0, \tpsi+ \delta_0]$.

\bs

A result similar to Theorem 1 can be stated e.g. for Axiom A diffeomorphisms on basic sets. Recall that 
if  $F: M \longrightarrow M$ is a $C^1$  Axiom A  diffeomorphism  on a Riemannian manifold $M$, a non-empty subset 
$\Lambda$ of $M$ is called a {\it basic set} for $F$ if  $\Lambda$ is a locally maximal compact $F$-invariant subset of $M$ 
which is not a single orbit, $F$ is hyperbolic  and transitive on $\Lambda$, and the periodic points of $F$ in $\Lambda$ 
are dense in $\Lambda$ (see e.g. \cite{B} or Appendix III in \cite{PP}). It follows from the existence of Markov
partitions that there exists a two-sided subshift of finite type $\sigma : \sa \longrightarrow \sa$ and a continuous 
surjective map $\pi : \sa \longrightarrow \Lambda$ such that: (i) $F \circ \pi = \pi \circ \sigma$, and (ii) for every H\"older
continuous function $g$ on $\Lambda$, $f = g\circ \pi\in \ff_\theta$ for some $\theta \in (0,1)$ and $\pi$ is one-to-one
almost everywhere with respect to the equilibrium state of $f$. Given a H\"older continuous function $g$ on $\Lambda$, 
the rate function $I_g$ is naturally related to the rate function $I_f$ of $f = g\circ \pi$. On the other hand,
$f$ is cohomologous to a function $f' \in \ff_{\sqrt{\theta}}(\sa)$ which depends on forward coordinates only, so 
$f' \in \ff_{\sqrt{\theta}}(\saa)$. Applying Theorem 1 to $f'$ provides a similar result for $f$ and therefore for $g$.

For some hyperbolic systems, large deviation principles similar to (\ref{eq:1.1}), however with shrinking intervals,
have been established recently in \cite{PoS} and \cite{PeS}.

% Sect. 2
\section{Proof of Theorem 1}
\setcounter{equation}{0}

\subsection{The Ruelle-Perron-Frobenius Theorem}

For convenience of the reader we state here the part of the estimates in \cite{St} that will be used in this section.

\bs

\noindent
{\bf Theorem 2} (Ruelle-Perron-Frobenius)
{\it Let  the $s_0\times s_0$ matrix $A$ and $M > 0$ be as in Sect. 1,  let $f\in \ff_\theta(\saa)$ be real-valued, and let
$b_f = \max \{ 1, | f |_\theta\}$. Then:}

\ms

   (i) {\it There exist a unique $\lambda = \lambda_f > 0$, a probability measure $\nu = \nu_f$ on $\saa$ and a
positive function $h = h_f \in \ff_\theta(\saa)$ such that $L_f h = \lambda\, h$ and $\di \int h\, d\nu = 1$.
The spectral radius of $L_f$ as an operator on $\ff_\theta(\saa)$ is $\lambda$, and 
its essential spectral radius  is $\theta\, \lambda$. The eigenfunction $h$ satisfies
$$\|h\|_\theta \leq \frac{6\,s_0^{M}\, b_f}{\theta^{2}(1-\theta)}\, e^{4 \, b_f /(1-\theta)}\, e^{2M | f |_\infty}$$
and
$$\min \, h \geq \frac{1}{ e^{2\, \,b_f/(1-\theta)}\, s_0^{M }\, e^{2M | f |_\infty} \;.}$$
Moreover,
$$\frac{\min \, h}{|h|_\infty}\, \lambda^n \leq L_f^n 1 \leq \frac{|h|_\infty}{\min\, h} \, \lambda^n\;,$$
for any integer $n \geq 0$.}

\ms

 (ii) {\it  The probability measure $\hat{\nu} = h\,\nu$
(this is the so called {\it equilibrium state} of $f$) is $\sigma$-invariant and 
$\hat{\nu} = \nu_{\hf}$, where $\hf = f - \log (h\circ \sigma) + \log h - \log \lambda$. Moreover $L_{\hf} 1 = 1$.}

\ms

(iii) {\it For  every  $g\in \ff_\theta(\saa)$ and every integer $n \geq 0$ we have
$$\left\| \frac{1}{\lambda^n}\; L_f^n g - h\, \int g \, d\nu \right\|_\theta \leq  D \,  \rho^n\, \|g\|_\theta\;,$$
where we can take
$$\rho = \left( 1- \frac{1-\theta}{4\, 
s_0^{2M}\, e^{\frac{8\,\theta \, b_f }{1-\theta}}\, e^{4M\,  | f |_\infty} } \right)^{\frac{1}{2M}} \in (0,1) $$
and}
$$D =  10^8\, \frac{b_f^7}{\theta^{10}\,(1-\theta)^{8}}\,  s_0^{17 M}\,e^{40  \, b_f/(1-\theta)}\,e^{33  M\,| f |_\infty} . $$

\bs

\noindent
{\bf Remark.} The constants  that appear in the above estimates are not optimal. The proof of Theorem 2 in \cite{St} follows that in
Sect. 1.B in \cite{B} with a more careful analysis of the estimates involved. The main point here is that,
apart from their obvious dependence on parameters related to the subshift of finite type $\sigma : \saa \longrightarrow \saa$,
these constant can be taken to depend only on $| f |_\theta$ and $| f |_\infty$.

\bs

\subsection{Reductions}

Let $f \in \ff_\theta(\saa)$ be the fixed function from Sect. 1 and let $\mu = \mu_f$ as before.
It follows from the properties of pressure (see e.g. \cite{R} or \cite{PP}) that $\Pr(g + c) = \Pr(g) + c$ for every
continuous function $g$ and every constant $c \in \R$. Thus, replacing $f$ by $f - \Pr(f)$, we may assume that
$\Pr(f) = 0$. Moreover, if $g$ and $h$ are cohomologous continuous functions on $\saa$, then $\Pr(g) = \Pr(h)$
and the equilibrium states $\mu_g$ of $g$ and $\mu_h$ of $h$ on $\saa$ coincide. Since $f$
is cohomologous to a function $\phi \in \ff_\theta(\saa)$ with $L_\phi 1 =1$ (see e.g. \cite{PP}), it is enough to 
prove the main result with $f$ replaced by such a function $\phi$. Moreover, $|\phi|_\infty$ and $|\phi|_\theta$
can be estimated by means of $|f- \Pr(f)|_\infty$ and $|f|_\theta$ (see e.g. part (ii) of Theorem 2 above).

So, from now on we will assume that  $L_\phi 1 = 1$.
It then follows that $\Pr(\phi) = 0$. Let $\mu = \mu_\phi$ be the {\it equilibrium state} of $\phi$ on $\saa$.

For the proof of Theorem 1 we may assume that $\psi \geq 0$. Indeed, assuming the statement of the theorem is true in
this case, suppose $\psi$ takes negative values. Set $\psi_1 = \psi + c$, where $c = - \min \psi$. 
Then $\psi_1 \geq 0$. Moreover, $\tpsi_1 = \int \psi_1\, d\mu = \tpsi + c$, $B_{\psi_1} = B_\psi$, and for $p_1 = p+c$
we have 
$$\Gamma_1(q) = p_1 q - \Pr(\phi+ q\psi_1) = (p+c)q - \Pr(\phi + q\psi + qc) = pq - \Pr(\phi + q\psi) = \Gamma(q) $$
for all $q \in \R$. Therefore (\ref{eq:1.2}) implies
$$I_{\psi}(p) = \sup \{pq - \Pr(\phi + q\psi) : q\in \R \} = \sup \{p_1q - \Pr(\phi + q\psi_1) : q\in \R \} = I_{\psi_1}(p_1) .$$
Moreover, if $0 < \delta _0 < B_\psi = B_{\psi_1}$, then $p \notin [\tpsi - \delta_0, \tpsi+\delta_0]$ is equivalent to 
$p_1 = p+c \notin [\tpsi_1 - \delta_0, \tpsi_1+\delta_0]$. 
Since $|\psi_1|_\theta = |\psi|_\theta$ and $|\psi_1|_\infty \leq 2 |\psi|_\infty$, using Theorem 1 for $I_{\psi_1}(p_1)$
and changing appropriately the value of the constant $q_0$, we get a similar estimate for  $I_\psi(p)$.

% Sect. 2.2
\subsection{Proof of Theorem 1 for $\psi \geq 0$}

From now on we will assume that $\phi, \psi \in \ff_\theta (\saa)$ are fixed real-valued functions such that $\psi \geq 0$, 
$\psi$ is not cohomologous to a constant, and 
\begin{equation} \label{eq:2.1}
L_\phi 1 = 1 .
\end{equation}

Given any $q\in \R$, set
$$f_q = \phi + q \psi \quad, \quad L_q = L_{f_q} .$$
In what follows we will assume 
\begin{equation}  \label{eq:2.2}
|q| \leq  q_0 \leq \frac{1}{b} 
\end{equation}
for some constant $q_0 > 0$ which will be chosen below. 
Then $|f_q|_\theta \leq |\phi|_\theta + 1$ for all $q$ with (\ref{eq:2.2}), and also $|f_q|_\infty \leq |\phi|_\infty + |\psi|_\infty$. Thus,
setting
$$C_0 =  \|\phi\|_\theta +  2 \max \{  |\psi|_\infty , 1\} \geq 1,$$
we have 
\begin{equation}  \label{eq:2.3}
\|f_q\|_\theta \leq C_0 \quad, \quad |q| \in [0,q_0] .
\end{equation}

Let $\nu_q$ be the {\it probability measure on $\saa$ with}
\be  \label{eq:2.4}
L^*_q \nu_q = \lambda_q \nu_q ,
\ee
where $\lambda_q$ is the {\it maximal eigenvalue} of $L_q = L_{f_q}$, and let $h_q > 0$ be a {\it corresponding normalised eigenfunction},
i.e. $h_q \in \ff_\theta(\saa)$, $L_qh_q = \lambda_q h_q$ and $\di \int h_q \, d\nu_q = 1$. Then $\mu_q = h_q \nu_q$ is the 
{\it equilibrium state} of $f_q$, i.e. $\mu_q = \mu_{\phi + q \psi}$. Clearly $h_0 = 1$ and $\mu_0 = \mu$.

Using the uniform estimates in Theorem 2 above, it follows from (\ref{eq:2.3}) that there exist
constants $D \geq 1$ and $\rho \in (0,1)$, depending on $C_0$ but not on $q_0$, such that
\be  \label{eq:2.5}
\left\| \frac{1}{\lambda^n_q} L^n_q g - h_q \; \int g\, d\nu_q \right\|_\theta \leq D\, \rho^n\, \|g\|_\theta
\ee
for all integers $n \geq 0$, all functions $g \in \ff_\theta(\saa)$ and all $q$ with $|q| \in [0,q_0]$.

{\bf Set $L = L_\phi$.} Given $x \in \saa$ and $m \geq 0$, set
$g_m(x) = g(x) + g(\sigma x) + \ldots + g(\sigma^{m-1} x) .$

It follows from (\ref{eq:2.4}) with $g = 1$ that $\di \lambda_q= \int L_q 1 \, d\nu_q$. Now
$$(L_q 1)(x)  = \sum_{\sigma y = x} e^{f_q(y)}  = \sum_{\sigma y = x} e^{\phi(y)+q\psi(y)}
 \leq e^{q_0 |\psi|_\infty} \; (L \,1)(x) =  e^{q_0 |\psi|_\infty}  $$
for all $x\in \saa$ implies $\lambda_q \leq e^{q_0|\psi|_\infty}$. 
Similarly,  $\lambda_q \geq e^{-q_0|\psi|_\infty}$.  Thus,
\be  \label{eq:2.6}
e^{-q_0 C_0} \leq \lambda_q \leq e^{q_0 C_0}  \quad, \quad |q| \leq q_0 .
\ee

To estimate $h_q$ for $q$ with (\ref{eq:2.2}), first use (\ref{eq:2.5}) with $g =1$ to get
$$\left\| \frac{1}{\lambda^n_q} L^n_q 1 - h_q\right\|_\theta \leq D\, \rho^n  .$$
Using (\ref{eq:2.1}), this gives
\begin{eqnarray*}
h_q (x)
& \leq & \frac{(L^n_q1)(x)}{\lambda^n_q} + D \rho^n
= \frac{1}{\lambda^n_q} \, \sum_{\sigma^ny = x} e^{(\phi+ q\psi)_n(y) } + D\, \rho^n\\
& \leq & \frac{e^{q_0C_0 n}}{\lambda^n_q}\, (L\, 1)(x) + D\, \rho^n \leq  e^{2q_0 C_0 n} + D\, \rho^n \nonumber
\end{eqnarray*}
for all $x \in \saa$ and $n \geq 0$. Similarly, 
$$h_q \geq \frac{e^{-q_0C_0 n}}{\lambda^n_q}\, (L \, 1) -  D\, \rho^n \geq  e^{-2 q_0 C_0 n} -  D\, \rho^n$$
for all $n \geq 0$. Thus,
\be  \label{eq:2.7}
\max\{ 0  , e^{- 2q_0 C_0 n} - D \rho^n\} \leq h_q \leq e^{2q_0 C_0 n} + D \rho^n \quad , \quad n \geq 0 \;,\; |q| \leq q_0 .
\ee

From now on we will assume that {\bf $p \notin [\tpsi - \delta_0, \tpsi + \delta_0]$ is fixed.}
Consider the function
$$\Gamma(q) = pq - \Pr(\phi + q \psi) \quad , \quad q \in \R .$$
Then
$I(p) = \sup_{q\in \R} \Gamma(q) .$

Clearly, $\Gamma(0) = 0$ and moreover by (\ref{eq:1.3}),
\be  \label{eq:2.8}
\Gamma'(q) = p - \int \psi\; d\mu_{\phi + q\psi}  .
\ee
In particular, $\Gamma'(0) = p - \tpsi \notin [-\delta_0, \delta_0]$.
 
We will now estimate the integral in the  right-hand-side of (\ref{eq:2.8}).
Let $\alpha > 0$ be the constant so that $\rho_1 = \max \{ \rho, \theta \} = e^{-\alpha}$.

\bs

\noindent
{\bf Lemma 1.}
{\it Assume that $\psi \geq 0$ on $\saa$ and $0 < \delta_0 <  B_\psi = \tpsi$. Set 
\be  \label{eq:2.9}
q_0 =  \min \left\{ \frac{\delta_0}{100 C_0^2 n_0} , \frac{1}{b} \right\} ,
\ee
where $n_0$ is the integer with
\be  \label{eq:2.10}
n_0  -1 \leq \frac{1}{\alpha} \left| \log \frac{\delta_0}{16C_0D} \right| <  n_0 .
\ee
Then  $\Gamma(q_0) \geq \frac{\delta_0 q_0}{2}$ and $\Gamma(-q_0) \geq \frac{\delta_0 q_0}{2}$.}

\bs

\noindent
{\bf Proof.}
For any $q \in [0,q_0]$ and any integer $n \geq 0$, (\ref{eq:2.4}), (\ref{eq:2.6}) and (\ref{eq:2.7}) yield
\begin{eqnarray*}
\int \psi\, d\mu_q 
& =     & \int \psi  h_q\; d\nu_q = \frac{1}{\lambda^n_q} \int L^n_q(\psi h_q)\, d\nu_q\\
& =     & \frac{1}{\lambda^n_q} \int \sum_{\sigma^n y = x} e^{(\phi+ q\psi)_n(y)} \psi(y) h_q(y) \, d\nu_q(x)\\
& \leq & e^{2q_0C_0 n}\, (e^{2qC_0 n} + D \rho^n)\; \int L^n \psi\, d\nu_q .
\end{eqnarray*}
It follows from (2.5) with $q = 0$ and $g = \psi$ and the choice of $C_0$ that
\be  \label{eq:2.11}
\left| L^n \psi - \int \psi \, d\mu \right| = \left| L^n \psi - \int \psi \, d\nu \right| \leq D\, \rho^n C_0 ,
\ee
therefore $L^n\psi \leq \tpsi + C_0 D \, \rho^n$. Combining this with the above gives
\be  \label{eq:2.12}
\int \psi\, d\mu_q \leq e^{2q_0C_0 n}\, (e^{2q_0C_0 n} + D \rho^n) (\tpsi + C_0 D\, \rho^n) .
\ee

Let $n_0 = n_0(f,\theta, \delta_0) \geq 1$ be {\it the integer such that}
\be  \label{eq:2.13}
e^{-n_0 \alpha } < \frac{\delta_0}{16 C_0 D} \leq e^{- (n_0-1)\alpha} .
\ee
Then $- n_0  \alpha <  \log \frac{\delta_0}{16 C_0D} \leq -(n_0-1)  \alpha$, so $n_0$ satisfies (\ref{eq:2.10}).
With this choice of $n_0$ define $q_0$ by (\ref{eq:2.9}).  Then  for $q\in [0,q_0]$ we have
$12 qC^2_0n_0 \leq \delta_0/8$ and so $12 qC_0 n_0 \leq 1$. It now follows from (\ref{eq:2.12}) with $q\in [0,q_0]$
and $n = n_0$, $0 < \delta_0 \leq B_\psi = \tpsi \leq C_0$, (\ref{eq:2.13}) and the fact that $e^{x} \leq 1+ 3x$ for $x\in [0,1]$ that
\begin{eqnarray*}
\int \psi\, d\mu_q 
& \leq & (e^{4q_0C_0 n_0} + D e^{2q_0C_0 n_0} e^{-\alpha n_0}) (\tpsi + C_0 D\, e^{-\alpha n_0})\\
& \leq & \left(1 + 12 q_0C_0 n_0 + (1+ 6 q_0C_0n_0)\frac{\delta_0}{16C_0}\right) \left(\tpsi + \frac{\delta_0}{16}\right)\\
& \leq & \tpsi + 12 q_0C^2_0 n_0 + (1+ 6 q_0C_0n_0)\frac{\delta_0}{16} + \left(2 + 2 \frac{\delta_0}{16C_0}\right) \frac{\delta_0}{16}\\
& \leq & \tpsi + \frac{\delta_0}{8} + \frac{\delta_0}{8} + \frac{3\delta_0}{16} < \tpsi + \frac{\delta_0}{2} .
\end{eqnarray*}

Thus, in the case $p \geq \tpsi + \delta_0$, it follows from (\ref{eq:2.8}) that $\Gamma'(q) \geq \frac{\delta_0}{2}$ for all $q\in [0,q_0]$, 
and therefore $\Gamma(q_0) \geq \frac{\delta_0 q_0}{2} .$

Next, assume that $p \leq \tpsi - \delta_0$. We will now estimate $\int \psi \, d\mu_q$ from below for $q \in [-q_0,0]$. 
As in the previous estimate, using (\ref{eq:2.6}) and (\ref{eq:2.7}), for such $q$ we get
\begin{eqnarray*}
\int \psi\, d\mu_q 
& =     & \int \psi  h_q\; d\nu_q = \frac{1}{\lambda^{n_0}_q} \int L^{n_0}_q(\psi h_q)\, d\nu_q\\
& =     & \frac{1}{\lambda^{n_0}_q} \int \sum_{\sigma^{n_0} y = x} e^{(\phi+ q\psi)_{n_0}(y)} \psi(y) h_q(y) \, d\nu_q(x)\\
& \geq & e^{-2q_0C_0 n_0}\, (e^{-2q_0C_0 n_0} - D \rho^{n_0})\; \int L^{n_0} \psi\, d\nu_q .
\end{eqnarray*}
Notice that by the choice of $q_0$ and $n_0$ we have $e^{-2q_0C_0 n_0} - D \rho^{n_0} > 0$. In fact, it follows
from $e^{-x} > 1-x$ for $x > 0$ that $e^{-2q_0C_0 n_0} > 1 - 2q_0C_0 n_0$, while (\ref{eq:2.13}) implies 
$D\rho^{n_0} < \frac{\delta_0}{16 C_0}$. Thus,
$$e^{-2q_0 C_0 n_0} - D\, \rho^{n_0} > 1 - 2q_0 C_0 n_0 - \frac{\delta_0}{16 C_0} > 1- \frac{\delta_0}{8 C_0}.$$
On the other hand, (\ref{eq:2.11}) yields $\int L^{n_0} \psi \, d\nu_q \geq \tpsi - D C_0 \rho^{n_0} > \tpsi - \frac{\delta_0}{16}$.
Hence for $ q \in [-q_0,0]$ we get
\begin{eqnarray*}
\int \psi\, d\mu_q 
& \geq & (1 - 2q_0C_0 n_0) \left(1- \frac{\delta_0}{8 C_0}\right) \left(\tpsi - \frac{\delta_0}{16}\right)\\
& \geq & \left(1 - 2 q_0C_0 n_0 - \frac{\delta_0}{8C_0}\right) \left(\tpsi - \frac{\delta_0}{16}\right)\\
& \geq & \tpsi  - \tpsi \left(2 q_0C_0n_0 + \frac{\delta_0}{8C_0} \right) - \frac{\delta_0}{16}
 \geq  \tpsi - \frac{\delta_0}{50} - \frac{\delta_0}{8} - \frac{\delta_0}{16} >  \tpsi - \frac{\delta_0}{2} .
\end{eqnarray*}
Thus, for $q\in [-q_0,0]$ we have 
$$\Gamma'(q) = p - \int \psi \, d\mu_q \leq \tpsi - \delta_0 - (\tpsi -\delta_0/2) \leq - \frac{\delta_0}{2} ,$$
and therefore $\Gamma(-q_0) \geq \frac{\delta_0 q_0}{2}$.
\endofproof

\bs

\noindent
{\it Proof of Theorem} 1. Assume again that $\psi \geq 0$.  Let $p \geq \tpsi + \delta_0$. Then $I(p) = \sup_{q\in \R} \Gamma(q)$, so by Lemma 1,
$I(p) \geq \Gamma(q_0) \geq \frac{\delta_0 q_0}{2}$. Similarly, for $p \leq \tpsi - \delta_0$ we get $I(p) \geq \frac{\delta_0 q_0}{2}$.

As explained in Sect. 2.2, the case of an arbitrary real-valued $\psi\in \ff_\theta(\saa)$ follows from the case $\psi \geq 0$.
\endofproof

\bs

\noindent
{\bf Acknowledgement:} Thanks are due to the referees for their valuable comments and\\ suggestions.

\bs

\end{document}